\documentclass[11pt]{amsart}

\usepackage{amssymb}
\usepackage{enumitem}
\usepackage{graphicx}
\usepackage{verbatim}
\usepackage{epstopdf}
\usepackage{color}

\newcommand{\eps}{\varepsilon}
\newcommand{\ph}{\varphi}
\newcommand{\NN}{\mathbb{N}}
\newcommand{\ZZ}{\mathbb{Z}}
\newcommand{\RR}{\mathbb{R}}

\newcommand{\MMM}{\mathcal{M}}

\newcommand{\Mac}{\MMM^{\mathrm{ac}}}

\newcommand{\KK}{\mathbf{K}}
\newcommand{\II}{\mathbf{I}}
\newcommand{\JJ}{\mathbf{J}}

\newcommand{\PPP}{\mathcal{P}}
\newcommand{\QQQ}{\mathcal{Q}}
\newcommand{\RRR}{\mathcal{R}}

\newcommand{\ulim}{\varlimsup}

\newcommand{\bl}{{\overline{\lambda}}}

\DeclareMathOperator{\graph}{graph}

\theoremstyle{plain}
\newtheorem{theorem}{Theorem}[section]

\newtheorem{conjecture}[theorem]{Conjecture}

\theoremstyle{definition}
\newtheorem{definition}[theorem]{Definition}

\theoremstyle{remark}

\newtheorem{remark}[theorem]{Remark}

\numberwithin{equation}{section}

\title[The geometric approach for constructing SRB measures]{The geometric approach for constructing Sinai-Ruelle-Bowen measures}
\author{Vaughn Climenhaga}
\address{Department of Mathematics \\ University of Houston \\ Houston, TX  77204, USA}
\email{climenha@math.uh.edu}
\urladdr{http://www.math.uh.edu/$\sim$climenha/}

\author{Stefano Luzzatto}
\address{Abdus Salam International Centre for Theoretical Physics (ICTP), Strada Costiera 11, Trieste, Italy}
\email{luzzatto@ictp.it}

\author{Yakov Pesin}
\address{Department of Mathematics \\ McAllister Building \\ Pennsylvania State University \\ University Park, PA 16802, USA}
\email{pesin@math.psu.edu}
\urladdr{http://www.math.psu.edu/pesin/}

\begin{document}

\date{\today}

\begin{abstract}
An important class of `physically relevant' measures for dynamical systems with hyperbolic behavior is given by Sinai--Ruelle--Bowen (SRB) measures. We survey various techniques for constructing SRB measures and studying their properties, paying special attention to the geometric `push-forward' approach.  After describing this approach in the uniformly hyperbolic setting, we review recent work that extends it to non-uniformly hyperbolic systems.
\end{abstract}

\thanks{V.C.\ was partially supported by NSF grant DMS-1362838.  Ya.P. was partially supported by NSF grant DMS--1400027. V.C. and Ya.P. would like to thank Erwin Schr\"odinger Institute and ICERM where the part of the work was done for their hospitality.} 

\maketitle

\emph{It is our great pleasure to have this survey included in this special issue dedicated to the 80th birthday of the great dynamicists D. Ruelle and Y. Sinai. We take this opportunity to acknowledge the tremendous impact that their work has had and continues to have in this field of research.} 

\section{Introduction}

Let $f\colon M\to M$ be a $C^{1+\alpha}$ diffeomorphism of a compact smooth Riemannian manifold $M$, and $U\subset M$ an open subset with the property that 
$\overline{f(U)}\subset U$. Such a set $U$ is called a \emph{trapping region} and the set
$\Lambda=\bigcap_{n\ge 0} f^n(U)$ a \emph{topological attractor} for $f$. We allow the case $\Lambda=M$. It is easy to see that $\Lambda$ is compact, $f$-invariant, and maximal (i.e., if $\Lambda'\subset U$ is invariant, then $\Lambda'\subset\Lambda$). We want to study the statistical properties of the dynamics in \(U\). Let \(m\) denote normalized Lebesgue measure on \(  M  \) and let $\mu$ be an arbitrary probability  measure on 
$\Lambda$. The set 
$$
B_\mu=\biggl\{x\in U: \frac{1}{n}\sum_{k=0}^{n-1}h(f^k(x))\to\int_{\Lambda} h\,d\mu \text{ for any } h\in C^1(M)\biggr\}
$$ 
is called \emph{basin of attraction} of $\mu$. We say that $\mu$ is a \emph{physical measure} if $m(B_\mu)>0$. An attractor with a physical measure is often referred to as a \emph{Milnor attractor}, see \cite{Mil, Kane}.

The simplest example of a physical measure is when \(  \Lambda=\{p\}  \) is a single fixed point, in which case the Dirac-delta measure \(  \delta_{p}  \) is a physical measure and \(  U\subseteq B_{\delta_{p}}  \). A less trivial case occurs when \(  \Lambda=M  \) and 
\(  \mu  \) is an ergodic invariant probability measure with  \(  \mu\ll m  \). Then the invariance and ergodicity imply, by Birkhoff's Ergodic Theorem, that \(  \mu(B_{\mu})=1  \) and thus the absolute continuity immediately implies that \(  m(B_{\mu})>0  \). 

Both of these cases are, however, quite special, and a more general situation is when 
\( \Lambda  \) is a non-trivial attractor and \( m(\Lambda)=0 \), in which case any invariant measure is necessarily singular with respect to Lebesgue. In the 1970's Sinai, Bowen, and Ruelle constructed a special kind of physical measures, which are now called \emph{SRB} or \emph{Sinai-Ruelle-Bowen} measures, to deal with precisely this situation in the special case in which \(  \Lambda  \) is \emph{uniformly hyperbolic}.  One can make sense of the definition of SRB measure in much more general cases, and a large amount of research has been devoted in the last several decades to establishing the existence of SRB measure for attractors \(  \Lambda  \) which are not necessarily uniformly hyperbolic. The purpose of this note is to survey the results  and techniques which have been used to prove the existence of SRB measures in increasingly general situations, with an emphasis on the geometric `push-forward' approach.

\medskip

\emph{Overview of the paper:}
In Section \ref{sec:definition} we give the definition of an SRB measures and state some of its basic properties, such as that of being a physical measure. In Section \ref{sec:approach} we briefly describe some of the main strategies which have been used to construct SRB measures in various cases. In the remaining sections we discuss in a little more details the applications of some of these strategies to various classes of attractors: uniformly hyperbolic attractors, as originally considered by Sinai, Ruelle, and Bowen, in Section \ref{hyp-attr}, partially hyperbolic attractors in Section \ref{sec:part-hyp}, attractors with dominated splittings in Section \ref{dom-split}, non-uniformly hyperbolic attractors in Section \ref{non-uni-attr} and, finally, uniformly hyperbolic attractors with singularities in Section \ref{att-sing}. 

\section{Definition of SRB measures}
\label{sec:definition}

\subsection{Hyperbolic measures} 
We recall some important facts from non-uniform hyperbolicity theory, referring the reader to \cite{BP} for more details. Given $x\in\Lambda$ and $v\in T_xM$, the \emph{Lyapunov exponent} of $v$ at $x$ is defined by
$$
\chi(x,v)=\limsup_{n\to\infty}\frac1n\log\|df^nv\|, \quad x\in M, v\in T_xM.
$$
The function $\chi(x,\cdot)$ takes on finitely many values, $\chi_1(x)\le\dots\le\chi_p(x)$, where $p=\dim M$. The values of the Lyapunov exponent are invariant functions, i.e., $\chi_i(f(x))=\chi_i(x)$ for every $i$.

A Borel invariant measure $\mu$ on $\Lambda$ is \emph{hyperbolic} if $\chi_i(x)\ne 0$ and $\chi_1(x)<0<\chi_p(x)$; that is 
$$
\chi_1(x)\le\dots\chi_k(x)<0<\chi_{k+1}(x)\le\dots\le\chi_p(x)
$$
for some $k(x)\ge 1$. If $\mu$ is ergodic, then $\chi_i(x)=\chi_i(\mu)$ for almost every 
$x$. The non-uniform hyperbolicity theory (see \cite{BP}) ensures that for a hyperbolic measure $\mu$ and almost every $x\in\Lambda$, the following are true. 
\begin{enumerate}
\item There is a splitting $T_xM=E^s(x)\oplus E^u(x)$ where 
$$
\begin{aligned}
E^s(x)&=E_f^s(x)=\{v\in T_xM: \chi(x,v)<0\}, \\
E^u(x)&=E^u_f(x)=E_{f^{-1}}^s(x)
\end{aligned}
$$ 
are \emph{stable} and \emph{unstable subspaces} at $x$; they satisfy  
\begin{enumerate}
\item $dfE^s(x)=E^s(f(x))$ and $dfE^u(x)=E^u(f(x))$;
\item $\angle(E^s(x),E^u(x))\ge K(x)$ for some Borel function $K(x)>0$ on $\Lambda$ that satisfies condition \eqref{CK} below.
\end{enumerate}
\item There are \emph{local stable} $V^s(x)$ and \emph{local unstable} $V^u(x)$ manifolds at $x$; they satisfy 
$$
\begin{aligned}
d(f^n(x),f^n(y))&\le C(x)\lambda^n(x)d(x,y), \,&y\in V^s(x), \,n\ge 0,\\
d(f^{-n}(x),f^{-n}(y))&\le C(x)\lambda^n(x)d(x,y), \, &y\in V^u(x), \,n\ge 0
\end{aligned}
$$
for some Borel function $C(x)>0$ on $\Lambda$ that satisfies \eqref{CK}, and some Borel $f$-invariant function $0<\lambda(x)<1$.
\item There are the \emph{global stable} $W^s(x)$ and \emph{global unstable} $W^u(x)$ manifolds at $x$ (tangent to $E^s(x)$ and $E^u(x)$, respectively) so that
$$
W^s(x)=\bigcup_{n\geq 0}f^{-n}(V^s(f^n x)), \quad W^u(x)=\bigcup_{n\geq 0}f^n(V^u(f^{-n}x));
$$
these manifolds are invariant under $f$, i.e., $f(W^s(x))=W^s(f(x))$ and $f(W^u(x))=W^u(f(x))$.
\item\label{CK} The functions $C(x)$ and $K(x)$ can be chosen to satisfy
$$
C(f^{\pm1}(x))\le C(x)e^{\varepsilon(x)}, \quad K(f^{\pm1}(x))\ge K(x)e^{-\varepsilon(x)},
$$
where $\eps(x)>0$ is an $f$-invariant Borel function.
\item The \emph{size} $r(x)$ of local manifolds satisfies 
$r(f^{\pm1}(x))\ge r(x)e^{-\varepsilon(x)}$.
\end{enumerate}
One can show that $W^u(x)\subset\Lambda$ for every $x\in\Lambda$ (for which the global unstable manifold is defined).

Since $\lambda(x)$ is invariant, it is constant $\mu$-a.e.\ when $\mu$ is ergodic, so from now on we assume that $\lambda(x) = \lambda$ is constant on $\Lambda$.

Given $\ell>1$, define \emph{regular set} of level $\ell$ by
$$
\Lambda_\ell=\Bigl\{x\in\Lambda: C(x)\le\ell, \, K(x)\ge\frac{1}{\ell}\Bigr\}.
$$
These sets satisfy: 
\begin{itemize}
\item $\Lambda_\ell\subset\Lambda_{\ell+1}$, 
$\bigcup_{\ell\ge 1}\,\Lambda_\ell=\Lambda$;
\item the subspaces $E^{s,u}(x)$ depend continuously on $x\in\Lambda_\ell$; in fact, the dependence is H\"older continuous:
$$
d_G(E^{s,u}(x),E^{s,u}(y))\le M_\ell d(x,y)^\alpha,
$$ 
where $d_G$ is the Grasmannian distance in $TM$;
\item the local manifolds $V^{s,u}(x)$ depend continuously on $x\in\Lambda_\ell$; in fact, the dependence is H\"older continuous: 
$$
d_{C^1}(V^{s,u}(x),V^{s,u}(y))\le L_\ell d(x,y)^\alpha;
$$
\item $r(x)\ge r_\ell > 0$ for all $x\in\Lambda_\ell$.
\end{itemize}

\subsection{SRB measures} 
We can choose $\ell$ such that $\mu(\Lambda_\ell)>0$. For $x\in\Lambda_\ell$ and a small $\delta_\ell>0$ set
$$
Q_\ell(x)=\bigcup_{y\in B(x,\delta_\ell)\cap\Lambda_\ell}\, V^u(y).
$$
Let $\xi_\ell$ be the partition of $Q_\ell(x)$ by $V^u(y)$, and let $V^s(x)$ be a local stable manifold that contains exactly one point from each $V^u(y)$ in $Q_\ell(x)$.  Then there are \emph{conditional measures} $\mu^u(y)$ on each $V^u(y)$, and a \emph{transverse measure} $\mu^s(x)$ on $V^s(x)$, such that for any $h\in L^1(\mu)$ supported on $Q_\ell(x)$, we have
\begin{equation}\label{eqn:conditional}
\int h\,d\mu = \int_{V^s(x)} \int_{V^u(y)} h \,d\mu^u(y) \,d\mu^s(x).
\end{equation}
See \cite[\S1.5]{CK} and references therein for further details. Let $m_{V^u(y)}$ denote the leaf volume on $V^u(y)$.
\begin{definition}
A measure $\mu$ on $\Lambda$ is called an \emph{SRB measure} if $\mu$ is hyperbolic and for every $\ell$ with $\mu(\Lambda_\ell)>0$, almost every $x\in\Lambda_\ell$ and almost every $y\in B(x,\delta_\ell)\cap\Lambda_\ell$, we have the measure $\mu^u(y)$ is absolutely continuous with respect to the measure  $m_{V^u(y)}$.  
\end{definition}

For $y\in\Lambda_\ell$, $z\in V^u(y)$ and $n>0$ set
$$
\rho^u_n(y,z)=\prod_{k=0}^{n-1}\frac{\text{Jac}(df|E^u(f^{-k}(z)))}{\text{Jac}(df|E^u(f^{-k}(y)))}.
$$
One can show that for every $y\in\Lambda_\ell$ and $z\in V^u(y)$ the following limit exists
\begin{equation}\label{eq:density}
\rho^u(y,z)=\lim_{n\to\infty}\rho^u_n(y,z)=\prod_{k=0}^{\infty}\frac{\text{Jac}(df|E^u(f^{-k}(z)))}{\text{Jac}(df|E^u(f^{-k}(y)))}
\end{equation}
and that $\rho^u(y,z)$ depends continuously on $y\in\Lambda_\ell$ and $z\in V^u(y)$.
\begin{theorem}[\cite{BP}, Theorems 9.3.4 and 9.3.6]
If $\mu$ is an SRB measure on $\Lambda$, then the density $d^u(x,\cdot)$ of the conditional measure $\mu^u(x)$ with respect to the leaf-volume $m_{V^u(x)}$ on $V^u(x)$ is given by $d^u(x,y)=\rho^u(x)^{-1}\rho^u(x,y)$ where 
$$
\rho^u(x)=\int_{V^u(x)}\rho^u(x,y)\,dm^u(x)(y)
$$
is the normalizing factor.
\end{theorem}
In particular, we conclude that the measures $\mu^u(x)$ and $m_{V^u(y)}$ must be equivalent.

The idea of describing an invariant measure by its conditional probabilities on the elements of a continuous partition goes back to the classical work of Kolmogorov and especially later work of Dobrushin on random fields (see \cite{Dob}). Relation \eqref{eq:density} can be viewed as an analog of the famous Dobrushin-Lanford Ruelle equation in statistical physics, see \cite{LanRu} and \cite{Sin2}.

\subsection{Ergodic properties of SRB measures} 

Using results of nonuniform hyperbolicity theory one can obtain a sufficiently complete description of ergodic properties of SRB measures.  
\begin{theorem}\label{srb-ergodic}
Let $f$ be a $C^{1+\alpha}$ diffeomorphism of a compact smooth manifold $M$ with an attractor $\Lambda$ and let $\mu$ be an SRB measure on
 $\Lambda$. Then there are $\Lambda_0,\Lambda_1,\Lambda_2,\dots \subset \Lambda$ such that
\begin{enumerate}
\item $\Lambda=\bigcup_{i\ge 0}\Lambda_i$, \,\,
$\Lambda_i\cap\Lambda_j=\emptyset$;
\item $\mu(\Lambda_0)=0$ and $\mu(\Lambda_i)>0$ for $i>0$;
\item $f|\Lambda_i$ is ergodic for $i>0$;
\item for each $i>0$ there is $n_i>0$ such that
$\Lambda_i=\bigcup_{j=1}^{n_i}\Lambda_{i,j}$ where the union is disjoint (modulo $\mu$-null sets), $f(\Lambda_{i,j})=\Lambda_{i,j+1}$, 
$f(\Lambda_{n_i,1})=\Lambda_{i,1}$ and $f^{n_i}|\Lambda_{i,1}$ is Bernoulli;
\item if $\mu$ is ergodic, then the basin of attraction $B_\mu$ has positive Lebesgue measure in $U$.
\end{enumerate}
\end{theorem}
For smooth measures this theorem was proved by Pesin in \cite{Pes1} and an extension to the general case was given by Ledrappier in \cite{Led} (see also \cite{BP}). 

\begin{figure}[htbp]
\includegraphics{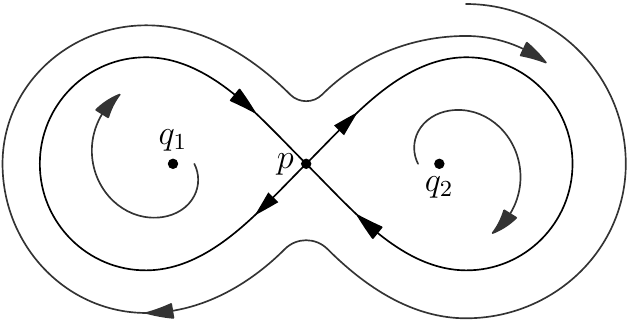}
\caption{A physical measure that is not SRB.}
\label{fig-eight}
\end{figure}

We stress that the final item of Theorem \ref{srb-ergodic} can be paraphrased as follows: \emph{ergodic SRB measures are physical}. The example of an attracting fixed point illustrates that the converse is not true; a more subtle example is given by the time-1 map of the flow illustrated in Figure \ref{fig-eight}, where the Dirac measure at the hyperbolic fixed point $p$ is a hyperbolic physical measure whose basin of attraction includes all points except $q_1$ and $q_2$. (In fact, by slowing down the flow near $p$, one can adapt this example so that $p$ is an indifferent fixed point and hence, the physical measure is not even hyperbolic.)

Returning to our discussion of SRB measures, one can show that a measure $\mu$ on 
$\Lambda$ of positive entropy is an SRB measure if and only if the entropy $h_\mu(f)$ of 
$\mu$ is given by the \emph{entropy formula}:
$$
h_\mu(f)=\int_\Lambda\,\sum_{\chi_i(x)>0}\,\chi_i(x)\,d\mu(x).
$$
For smooth measures (which are a particular case of SRB measures) the entropy formula was proved by Pesin \cite{Pes1} (see also \cite{BP}) and its extension to SRB measures was given by Ledrappier and Strelcyn \cite{LS}. The fact that a hyperbolic measure satisfying the entropy formula is an SRB measure was shown by Ledrappier \cite{Led}. \footnote{Our definition of SRB measure includes the requirement that the measure is hyperbolic. In fact, one can extend the notion of SRB measures to those that are \emph{non-uniformly partially hyperbolic}. In this case some Lyapunov exponents can be zero but there must be at least one positive Lyapunov exponent. It was proved by Ledrappier and Young \cite{LY} that a non-uniformly partially hyperbolic measure satisfies the entropy formula if and only if it is an SRB measure in this more general sense. We stress that a non-uniformly partially hyperbolic SRB measure may not be physical, i.e., its basin may be of zero Lebesgue measure.}

It follows from Theorem \ref{srb-ergodic} that $f$ admits at most countably many ergodic SRB measures. It is shown in \cite{HHTU} that a topologically transitive $C^{1+\alpha}$ surface diffeomorphism can have at most one SRB measure but the result is not true in dimension higher than two, see Section \ref{srb:uniq}.

\section{Approaches to the construction of SRB measures}\label{sec:approach} 

There exist at least three distinct arguments for the construction of SRB measures.  In this section we briefly describe these approaches and mention the different settings in which they have been applied.  Then beginning in \S\ref{hyp-attr} we give some more details of specific results and of the inner workings of each of the arguments.

The first approach, used by Sinai, Ruelle, and Bowen in their pioneering work, is based on the non-trivial fact that uniformly hyperbolic attractors admit a \emph{finite Markov partition} and consequently a \emph{symbolic coding} by a \emph{subshift of finite type (SFT)}. 
This coding makes it possible to translate questions about invariant measures for the diffeomorphism $f$ into the language of symbolic dynamics, and hence to borrow results from statistical mechanics regarding \emph{Gibbs measures} and \emph{equilibrium states} for certain \emph{potential functions} $\ph\colon \Lambda \to \RR$.  Of particular importance is the potential function \( \varphi = -\log |\det df^{u}| \), where \( \det df^{u} \) 
denotes the determinant of the differential of \( f \) restricted to the ``unstable'' subspace of the uniformly hyperbolic system; Gibbs measures for $\ph$ correspond to SRB measures for $(\Lambda,f)$.

We give a rough outline of the symbolic approach.
\begin{enumerate}
\item\label{Markov} Use a finite Markov partition to code $(\Lambda,f)$ by a two-sided SFT $\Sigma \subset A^\ZZ$, where $A$ is a finite alphabet.
\item\label{quotient} Pass from $\Sigma$ to the corresponding one-sided shift $\Sigma^+ \subset A^\NN$; roughly speaking, this corresponds to identifying points on $\Lambda$ that lie on the same local stable manifold.
\item\label{RPF} Consider a certain \emph{transfer operator} associated to $\Sigma^+$ and the potential function, and obtain a Gibbs measure in terms of the eigendata of this operator using Perron--Frobenius theory.
\item\label{Gibbs-SRB} Project this Gibbs measure on $\Sigma^+$ to an SRB measure on $\Lambda$.
\end{enumerate}
The extension of this argument to more general attractors which are not uniformly hyperbolic, or even to uniformly hyperbolic systems with singularities, is made more challenging by the fact that one cannot hope to have finite Markov partitions in more general settings. Even though one can construct \emph{countable} Markov partitions in some settings, the theory of Gibbs measures for shift maps on countable symbolic spaces is not as complete as for finite symbolic spaces, and thus some new ideas are needed in order to generalize this approach to the construction of SRB measures.

It is worth mentioning that the heart of the symbolic approach lies in the application of Perron--Frobenius theory by finding the appropriate Banach space on which the transfer operator acts with a spectral gap, and that for uniformly hyperbolic systems, this functional analytic strategy can in fact be carried out without relying on symbolic dynamics \cite{BGL}.  The key is to identify the right Banach space; roughly speaking one should consider objects that behave like smooth functions along the unstable direction, and like measures (or more generally, distributions) along the stable direction.

An alternative to the functional analytic approach, which is more ``geometric'', was developed in \cite{PS} to deal with partially hyperbolic attractors for which Markov partitions do not exist and the above symbolic approach fails (we describe the result in \cite{PS} more precisely in Section \ref{sec:part-hyp} below). The idea here is to follow the classical Bogolyubov-Krylov procedure for constructing invariant measures by pushing forward and a given reference measure. In our case the natural choice of a reference measure is the Riemannian volume $m$ restricted to the neighborhood $U$, which we denote by $m_U$. We then consider the sequence of probability measures 
\begin{equation}\label{eq:sequence1}
\mu_n=\frac{1}{n}\sum_{k=0}^{n-1}f_*^km_U.
\end{equation}
Any weak* limit point of this sequence of measures is called a \emph{natural measure} and while in general, it may be a trivial measure, under some additional hyperbolicity requirements on the attractor one obtains an SRB measure. 

For attractors with some hyperbolicity one can use a somewhat different approach which exploits the fact that SRB measures are absolutely continuos along unstable manifolds. To this end consider a point  $x\in\Lambda$, its local unstable manifold $V^u(x)$, and take the leaf-volume $m_{V^u(x)}$ as the reference measure for the above construction. Thus one studies the sequence of probability measures given by
\begin{equation}\label{eq:sequence2}
\nu_n=\frac{1}{n}\sum_{k=0}^{n-1}f_*^km_{V^u(x)}.
\end{equation}
The measures \( \nu_{n} \) are spread out over increasingly long pieces of the global unstable manifold of the point $f^n(x)$ and, in some situations, control over the geometry of the unstable manifold makes it possible to draw conclusions about the measures \( \nu_{n} \) by keeping track of their densities along unstable manifolds, and ultimately to demonstrate that any weak-star limit of \( \nu_{n} \) is in fact an SRB measure.  This approach applies to the uniformly hyperbolic setting, as an alternative to the symbolic coding approach, as we will discuss in Section \ref{hyp-attr}. It can also be applied to
significantly more general situations such as partially hyperbolic and non-uniformly hyperbolic settings, which we discuss in Sections \ref{sec:part-hyp}--\ref{non-uni-attr}.

Before discussing the key ideas necessary to extend either approach to the non-uniformly hyperbolic situation, we pause to describe the relationship between these two approaches.  Because the definition of SRB measure is so closely tied to the unstable direction, any approach to constructing SRB measures must somehow work with the unstable direction and ignore the stable one (at least at certain stages).
\begin{itemize}
\item In the symbolic approach, this is done in Step \ref{quotient} by passing from $\Sigma$ to $\Sigma^+$; for shift spaces, the location of a point along a stable manifold is encoded in the negative indices, while a location along an unstable manifold is encoded in the positive indices.
\item In the geometric approach, we privilege the unstable direction by working with $m_{V^u(x)}$ (instead of $m_{V^s(x)}$) and by keeping track of densities  along unstable manifolds.
\end{itemize}
Next one must take the dynamics of $f$ into account; after all, we are looking for an invariant measure.
\begin{itemize}
\item In the symbolic approach, this is done via the \emph{transfer operator}, which acts on the space of H\"older continuous functions on $\Sigma^+$, viewed as densities of a measure.
\item In the geometric approach, this is done via the time-averaging process in \eqref{eq:sequence2}.
\end{itemize}

Now suppose we wish to extend these approaches to settings with non-uniformly hyperbolic dynamics.  Both approaches rely on having uniform expansion and contraction properties, and as we will see in Sections \ref{sec:part-hyp}--\ref{non-uni-attr}, the key to extending the second (geometric) approach to the non-uniformly hyperbolic setting is to restrict one's attention to orbit segments where expansion and contraction occur uniformly (in a sense that will be made precise).  That is, instead of considering \emph{all} iterates $f^k(x)$, one considers only \emph{hyperbolic times}; that is, values of $k$ such that $f^k$ has some uniformly hyperbolic properties at $x$ (the set of such times depends on the point where the orbit starts).  

In a number of situations, it has proven possible to combine this idea with the symbolic coding techniques from the first approach, which leads us to the third approach, based on the concept of \emph{inducing}. If \( \Gamma\subseteq \Lambda \) is some appropriately chosen subset, and \( \tau: \Gamma \to \mathbb N \) is an \emph{inducing time}, or \emph{return time}, function, i.e.,  \( f^{\tau(x)}(x)\in \Gamma \) for every (or almost every) \( x\in \Gamma \), then we can define the induced map \( F: \Gamma \to \Gamma \) by \( F(x) = f^{\tau(x)}(x) \). The point here is that it may be possible to choose \( \Gamma \) with a much more amenable and regular geometric structure than the attractor \( \Lambda \) as a whole, and that it may be possible to choose a return time function \( \tau \) such that the induced map \( F \) has some good properties, for example \( F \) may be  (piecewise) uniformly hyperbolic. In these conditions we can then construct an SRB measures for \( F \) and by elementary and standard arguments, under some integrability conditions on the inducing time \( \tau \), use this to obtain an SRB measure for \( f \). The general notion of inducing is quite classical in ergodic theory but the specific application to SRB measures was first applied in the setting of certain H\'enon maps \cite{BY1} and developed as a general theory by Young \cite{You98}, from which specific kinds of induced maps required for SRB measures are often refereed to as \emph{Young towers}. 

Although the Young Tower approach is in principle more involved and complicated than the geometric ``push-forward'' approach, it has the advantage that the symbolic structure allows an application of techniques from functional analysis and spectra theory that give significantly more information about the structure and properties of SRB measure, including information on the statistical properties such as decay of correlations. 

Our main focus in the remainder of this paper will be the geometric approach.  The symbolic approach (for finite alphabets) is very well described in the original literature (see especially \cite{Bo2}). We will also discuss some results related to Young towers.

To carry out the construction of SRB measures using the push-forward geometric approach one needs certain information on the dynamics and geometry of ``unstable'' admissible manifolds and their images.\footnote{In the setting of non-uniform hyperbolicity, ``unstable'' admissible manifolds are used as substitutions for local unstable manifolds.} In particular, this includes \emph{hyperbolicity}. If it is uniform one can carry out the construction without too much trouble, see Section \ref{hyp-attr}. If hyperbolicity is not uniform, one may still hope to have the following.
\begin{enumerate}
\item Domination: if one of the directions does not behave hyperbolically, then it at least is still dominated by the other direction.
\item Separation: the stable and unstable directions do not get too close to each other, more precisely, there is a ``good'' way to control how close they can be. 
\end{enumerate}
In the case of attractors with dominated splittings (see Section \ref{dom-split}) these two conditions hold uniformly and so one only needs to control the asymptotic hyperbolicity (expansion and contraction along stable and unstable directions). In the more general non-uniformly hyperbolic case both domination and separation may fail at some points, and in order to control the geometry and dynamics of images of admissible manifolds, one needs to replace ``hyperbolicity'' with ``effective hyperbolicity'', see Section \ref{non-uni-attr}. 

\section{SRB measures for uniformly hyperbolic attractors}\label{hyp-attr}

\subsection{Definition of hyperbolic attractors} Consider a topological attractor $\Lambda$ for a diffeomorphism $f$ of a compact smooth manifold $M$. It is called \emph{(uniformly) hyperbolic} if for each $x\in\Lambda$ there is a decomposition of the tangent space $T_xM=E^s(x)\oplus E^u(x)$ and constants $c>0$, $\lambda\in(0,1)$ such that for each $x\in\Lambda$:
\begin{enumerate}
\item $\|d_xf^nv\|\le c\lambda^n\|v\|$ for $v\in E^s(x)$
and $n\ge0$;
\item $\|d_xf^{-n}v\|\le c\lambda^n\|v\|$ for $v\in E^u(x)$ and $n\ge0$.
\end{enumerate}
$E^s(x)$ and $E^u(x)$ are \emph{stable} and \emph{unstable subspaces} at $x$. One can show that $E^s(x)$ and $E^u(x)$ depend continuously on $x$. 

In particular, $\angle(E^s(x),E^u(x))$ is uniformly away from zero. In fact, $E^s(x)$ and $E^u(x)$ depend H\"older continuously on $x$.

For each $x\in\Lambda$ there are $V^s(x)$ and $V^u(x)$ \emph{stable} and \emph{unstable local manifolds} at $x$. They have uniform size $r$, depend continuously on $x$ in the $C^1$ topology and $V^u(x)\subset\Lambda$ for any 
$x\in\Lambda$. 

We describe an example of a hyperbolic attractor. Consider the solid torus 
$P=D^2\times S^1$. We use coordinates $(x,y,\theta)$ on $P$; $x$ and $y$ give the coordinates on the disc, and $\theta$ is the angular coordinate on the circle. Fixing parameters $a\in (0,1)$ and $\alpha,\beta\in (0,\min\{a,1-a\})$, define a map $f: P\to P$ by
$$
f(x,y,\theta) = (\alpha x + a\cos\theta, \beta y +a\sin\theta, 2\theta).
$$
$P$ is a trapping region and $\Lambda=\bigcap_{n\geq 0}f^n(P)$ is the attractor for $f$, known as the Smale-Williams solenoid.

\subsection{Existence of SRB measures for hyperbolic attractors}\label{exist-hyper-attr} The following result establishes existence and uniqueness of SRB measures for transitive hyperbolic attractors. 
\begin{theorem}\label{existence-srb}
Assume that $f$ is $C^{1+\alpha}$ and that $\Lambda$ is a uniformly hyperbolic attractor. The following statements hold:
\begin{enumerate}
\item Every limit measure of either the sequence of measures $\mu_n$ (given by \eqref{eq:sequence1}) or the sequence of measures $\nu_n$ (given by \eqref{eq:sequence2}) is an SRB measure on $\Lambda$.
\item There are at most finitely many ergodic SRB measures on $\Lambda$.
\item If $f|\Lambda$ is topologically transitive, then there is a unique SRB-measure $\mu$ on $\Lambda$, which is the limit of both sequences $\mu_n$ and $\nu_n$; moreover,
$B_\mu$ has full measure in $U$. 
\end{enumerate}
\end{theorem}
This theorem was proved by Sinai, \cite{Sin1} for the case of Anosov diffeomorphisms,   Bowen, \cite{Bo2} and Ruelle \cite{Ru} extended this result to hyperbolic attractors, and Bowen and Ruelle, \cite{BoRu} constructed SRB measures for Anosov flows. 

We outline a proof of this theorem to demonstrate how the geometric approach works; further details can be found in \cite{CDP}. Note that the geometric proof we give here is not the original proof given by Sinai, Bowen, and Ruelle, who used the symbolic approach. Given $x\in M$, a subspace $E(x)\subset T_x M$, and $a(x)>0$, the \emph{cone} at $x$ around $E(x)$ with angle $a(x)$ is
$$
K(x,E(x),a(x))=\{v\in T_x M\colon\measuredangle (v,E(x))<a(x)\}.
$$
There exists a neighborhood $\tilde{U}\subset U$ of the attractor $\Lambda$ and two continuous cone families $K^s(x)=K^s(x,E^s(x), a)$ and $K^u(x)=K^u(x,E^u(x), a)$ 
such that\footnote{Note that the subspaces $E^s(x)$ and $E^u(x)$ for $x\in\tilde U$ need not be invariant under $df$.}
\begin{align*}
\overline{df(K^u(x))}&\subset K^u(f(x)) \text{ for all }x\in\tilde{U},\\
\overline{df^{-1}(K^s(f(x)))}&\subset K^s(x) \text{ for all }x\in f(\tilde{U}).
\end{align*}
Let $W\subset U$ be an \emph{admissible manifold}; that is, a submanifold that is tangent to an unstable cone $K^u(x)$ at some point $x\in U$ and has a fixed size and uniformly bounded curvature. More precisely, fix constants $\gamma,\kappa,r>0$, and define a \emph{$(\gamma,\kappa)$-admissible manifold of size $r$} to be 
$V(x)=\exp_x \graph\psi$, where $\psi: B_{E^u(x)}(0,r)=B(0,r)\cap E^u(x)\to E^s(x)$ is $C^{1+\alpha}$ and satisfies
\begin{equation}\label{admissibe-man}
\begin{aligned}
\psi(0) &= 0 \text{ and } d\psi(0) = 0,\\
\|d\psi\| &:= \sup_{\|v\| < r} \|d\psi(v)\| \leq \gamma, \\
|d\psi|_\alpha &:=\sup_{\|v_1\|,\|v_2\|<r}\frac{\|d\psi(v_1)-d\psi(v_2)\|}{\|v_1-v_2\|^\alpha}\le\kappa.
\end{aligned}
\end{equation}
Write ${\bf I}=(\gamma,\kappa,r)$ for convenience and consider the space of admissible manifolds
\begin{multline*} 
\RRR_{\bf I}=\{\exp_x(\graph\psi)\colon x\in U,\psi\in C^1(B^u(0,r),E^s(x)) \text{ satisfies } \eqref{admissibe-man} \}.
\end{multline*}
Given an admissible manifold $W$, we consider a \emph{standard pair} $(W,\rho)$ where 
$\rho$ is a continuous ``density'' function on $W$. The idea of working with pairs of admissible manifolds and densities was introduced by Chernov and Dolgopyat \cite{CD} and is an important recent development in the study of SRB measures via geometric techniques. 

Now we fix $L>0$, write ${\bf K}=({\bf I},L)$, and consider the space of standard pairs
\[
\RRR'_{\bf K}=\{(W,\rho)\colon W\in\RRR_{\bf I}, \rho\in C^\alpha(W,[\tfrac 1L,L]), |\rho|_\alpha\le L\}.
\]
These spaces are compact in the natural product topology: the coordinates in 
$\RRR_{\bf I}$ are 
$$
\{x\in M, \psi\in C^1(B^u(0,r),E^s(x)) \text{ with } \|D\psi\|\le\gamma, \, 
|D\psi|_\alpha\le\kappa\}
$$
and the coordinates in $\RRR'_{\bf K}$ are
$$
\{x, \psi,\rho\in C^\alpha(W) \text{ with } \|\rho\|_\alpha\le L\}.
$$
A standard pair determines a measure $\Psi(W,\rho)$ on $\overline{U}$ in the obvious way:
$$
\Psi(W,\rho)(E) := \int_{E\cap W} \rho \,dm_W.
$$
Moreover, each measure $\eta$ on $\RRR'_{\bf K}$ determines a measure 
$\Phi(\eta)$ on $\overline{U}$ by
\begin{equation}\label{eqn:Phi2}
\begin{aligned}
\Phi(\eta)(E)&=\int_{\RRR'_{\bf K}}\Psi(W,\rho)(E)\,d\eta(W,\rho) \\
&=\int_{\RRR'_{\bf K}}\int_{E\cap W}\rho(x) \,d m_W(x) \,d\eta(W,\rho).
\end{aligned}
\end{equation}
(Compare this to \eqref{eqn:conditional} in the definition of conditional measures.)
Write $\mathcal{M}(\overline{U})$ and $\mathcal{M}(\RRR'_{\bf K})$ for the spaces of finite Borel measures on $\overline{U}$ and $\RRR'_{\bf K}$, respectively. It is not hard to show that $\Phi\colon\mathcal{M}(\RRR'_{\bf K})\to\mathcal{M}(\overline{U})$ is continuous; in particular, $\mathcal{M}_{\bf K}=\Phi(\mathcal{M}_{\le 1}(\RRR_{\bf K}'))$ is compact, where we write $\mathcal{M}_{\le 1}$ for the space of measures with total weight at most $1$. 

On a uniformly hyperbolic attractor, an invariant probability measure is an SRB measure if and only if it is in $\mathcal{M}_{\bf K}$ for some ${\bf K}$. 

Consider now the leaf volume $m_W$ on $W$ that we view as a measure on $\bar U$. Its evolution is the sequence of measures  
\begin{equation}\label{eqn:mun-mW}
\kappa_n=\frac1n\sum_{k=0}^{n-1}f_*^k m_W.
\end{equation}
By weak* compactness there is a subsequence $\kappa_{n_k}$ that converges to an invariant measure $\mu$ on $\Lambda$ which is an SRB measure.  

Consider the images $f^n(W)$ and observe that for each $n$, the measure $f_*^n m_W$ is absolutely continuous with respect to leaf volume on $f^n(W)$. For every $n$, the image $f^n(W)$ can be covered with uniformly bounded multiplicity (this requires a version of the Besicovitch covering lemma) by a finite number of admissible manifolds $W_i$, so that 
\begin{equation}\label{eqn:cvx-comb}
f_*^nm_W\text{ is a convex combination of measures }\rho_i\, dm_{W_i},
\end{equation}
where $\rho_i$ are H\"older continuous positive densities on $W_i$. 

We see from \eqref{eqn:cvx-comb} that $\mathcal{M}_{\bf K}$ is invariant under the action of $f_*$, and thus $\kappa_n\in\mathcal{M}_{\bf K}$ for every $n$. By compactness of 
$\mathcal{M}_{\bf K}$, one can pass to a subsequence $\kappa_{n_k}$ which converges to a measure $\mu\in\mathcal{M}_{\bf K}$, and this is the desired SRB measure.

Choosing $W=V^u(x)$, $x\in\Lambda$ we obtain that any limit measure of the sequence 
$\nu_n$ (see \eqref{eq:sequence2}) is an SRB measure. It is then not difficult to derive from here that any limit measure of the sequence $\mu_n$ (see \eqref{eq:sequence1}) is an SRB measure. 

In the particular case when $\Lambda=M$ (that is, $f$ is a $C^{1+\alpha}$ Anosov diffeomorphism) and $f$ is transitive, the above theorem guarantees existence and uniqueness of the SRB measure $\mu$ for $f$. Reversing the time we obtain the unique SRB measure $\nu$ for $f^{-1}$. One can show that $\mu=\nu$ if and only if $\mu$ is a smooth measure.

\section{SRB measures for partially hyperbolic attractors}
\label{sec:part-hyp}

\subsection{Definition of partially hyperbolic attractors} Consider a topological attractor 
$\Lambda$ for a diffeomorphism $f$ of a compact smooth manifold $M$. It is called \emph{(uniformly) partially hyperbolic} if for each $x\in\Lambda$ there is a decomposition of the tangent space $T_xM=E^s(x)\oplus E^c(x)\oplus E^u(x)$ and numbers 
$0<\lambda<\lambda_1\le\lambda_2<\lambda^{-1}$ and $c>0$ such that for $n\ge0$:
\begin{enumerate}
\item $\|d_xf^nv\|\le c\lambda^n\|v\|$ for $v\in E^s(x)$;
\item $c^{-1}\lambda_1^n\|d_xf^nv\|\le c\lambda_2^n\|v\|$;
\item $\|d_xf^{-n}v\|\le c\lambda^n\|v\|$ for $v\in E^u(x)$.
\end{enumerate}
Here $E^s(x)$, $E^c(x)$ and $E^u(x)$ are \emph{strongly stable}, \emph{central} and \emph{strongly unstable subspaces} at $x$. They depend (H\"older) continuously on $x$. In particular, the angle between any two of them is uniformly away from zero. 

For each $x\in\Lambda$ there are $V^s(x)$ and $V^u(x)$, the \emph{strongly stable} and \emph{strongly unstable local manifolds} at $x$. They have uniform size $r$, depend continuously on $x$ in the $C^1$ topology and $V^u(x)\subset\Lambda$ for any 
$x\in\Lambda$. 

A simple example of a partially hyperbolic attractor is a map which is the direct product of a map $f$ with a hyperbolic attractor $\Lambda$ and the identity map $\text{Id}$ of any manifold.  

\subsection{$u$-measures}\label{exist-hyper-attr1} 

In light of the absolute continuity condition for SRB measures, the following definition for a partially hyperbolic system is natural. A measure $\mu$ on $\Lambda$ is called a 
\emph{$u$-measure} if for every $x\in\Lambda$ and $y\in B(x,\delta)\cap\Lambda$, we have $\mu^u(y)\sim m_{V^u(y)}$. (Recall that $\mu\sim \nu$ if $\mu\ll \nu$ and $\nu\ll \mu$). 

\begin{theorem}[\cite{PS}]
Any limit measure of the sequence of measures $\mu_n$ (see \eqref{eq:sequence1}) is a $u$-measure and so is any limit measure of the sequence of measures $\nu_n$ (see \eqref{eq:sequence2}).  
\end{theorem}

This theorem is a generalization of Theorem~\ref{existence-srb} in the uniformly hyperbolic case and its proof uses the ``push-forward'' techniques. Indeed, recalling the definition of the push-forward of a measure we have 
\(f_{*}^{k}m_{V^{u}(x)}(A) = m_{V^{u}(x)}(f^{-k}(A))=m_{V^{u}(x)}(\{x: f^{k}(x)\in A\}) \) and hence, the measures \( f_{*}^{k}m_{V^{u}(x)}\) are supported on the image 
\( f^{k}(V^{u}(x)) \) of the starting chosen piece of local unstable manifold. If \( f \) is uniformly expanding along \( E^{u} \) one can divide up \( V^{u}(x) \) into pieces, each of which grows to large scale with bounded distortion at time \( k \), and thus 
\( f^{k}_{*}m_{V^{u}(x)} \) is supported on some collection of uniformly large unstable disks. Therefore, the same is true for measures \( \nu_{n} \) in \eqref{eq:sequence2} and this is the crucial property used to show that any limit measure \( \mu \) has absolutely continuous conditional measures along unstable leaves and therefore is an SRB measures (see Section \ref{exist-hyper-attr}).

One can prove the following basic properties for $u$-measures.
\begin{enumerate}
\item Any measure whose basin has positive volume is a $u$-measure, \cite{BDV}.
\item If there is a unique $u$-measure for $f$, then its basin has full volume in the topological basin of attraction, \cite{Dol}.
\item Every ergodic component of a $u$-measure is again a $u$-measure, \cite{BDV}. 
\end{enumerate}
The first of these says that for a partially hyperbolic attractor, every physical measure is a $u$-measure.  In particular, every SRB measure is a $u$-measure. What about the converse implications?  Are $u$-measures physical?  When is a $u$-measure an SRB measure? We address these questions next.

\subsection{$u$-measures with negative central exponents} We say that $f$ has \emph{negative (positive) central exponents} (with respect to $\mu$) if there exists an invariant subset $A\subset\Lambda$ with $\mu(A)>0$ such that the Lyapunov exponents $\chi(x,v)<0$ (respectively, $\chi(x,v)>0$) for every $x\in A$ and every vector 
$v\in E^c(x)$.

If $f$ has negative central exponents on a set $A$ of full measure with respect to a 
$u$-measure $\mu$, then $\mu$ is an SRB measure for $f$.
\begin{theorem}[\cite{BDPP}]\label{thm:neg}
Assume that $f$ has negative central exponents on an invariant set $A$ of positive measure with respect to a $u$-measure $\mu$ for $f$. Then the following statements hold:
\begin{enumerate}
\item Every ergodic component of $f|A$ of positive $\mu$-measure is open$\pmod 0$; in particular, the set $A$ is open$\pmod 0$ (that is there exists an open set $U$ such that 
$\mu(A\triangle U)=0$).
\item If for $\mu$-almost every $x$ the trajectory $\{f^n(x)\}$ is dense in $\text{supp}(\mu)$, then 
$f$ is ergodic with respect to~$\mu$.
\end{enumerate}
\end{theorem}
We provide the following criterion, which guarantees the density assumption in Statement (2) of the previous theorem. 
\begin{theorem}
Assume that for every $x\in\Lambda$ the orbit of the global strongly unstable manifold 
$W^u (x)$ is dense in $\Lambda$. Then for any $u$-measure $\mu$ on $\Lambda$ and 
$\mu$-almost every $x$ the trajectory $\{f^n(x)\}$ is dense in $\Lambda$.
\end{theorem}
This result is an immediate corollary of the following more general statement. Given 
$\varepsilon>0$, we say that a set is \emph{$\varepsilon$-dense} if its intersection with any ball of radius $\varepsilon$ is not empty.
\begin{theorem}[\cite{BDPP}]
Let $f$ be a $C^1$ diffeomorphism of a compact smooth Riemannian manifold $M$ possessing a partially hyperbolic attractor $\Lambda$. The following statements hold:
\begin{enumerate}
\item For every $\delta>0$ and every $\varepsilon\le\delta$ the following holds: assume that for every $x\in\Lambda$ the orbit of the global strongly unstable manifold $W^u(x)$ is 
$\varepsilon$-dense in $\Lambda$. Then for any $u$-measure $\mu$ on $\Lambda$ and 
$\mu$-almost every $x$ the trajectory $\{f^n(x)\}$ is $\delta$-dense in $\Lambda$.
\item Assume that for every $x\in\Lambda$ the orbit of the global strongly unstable manifold $W^u (x)$ is dense in $\Lambda$. Then $\text{supp}(\mu)=\Lambda$ for every $u$-measure 
$\mu$.
\end{enumerate}
\end{theorem}
In light of the `negative central exponent' hypothesis in Theorem \ref{thm:neg}, it is natural to ask whether a corresponding result holds for an attractor with positive central exponents. This case turns out to be more difficult since it is easier to handle non-uniformities in the contracting part of the dynamics than it is to handle non-uniformities in the expanding part of the dynamics. We discuss this situation in more detail in Section \ref{dom-split} but mention here that the study of $u$-measures with positive central exponents was carried out in \cite{ABV, ADLP} under the stronger assumption that there is a set of positive volume in a neighborhood of the attractor with positive central exponents.

\subsection{Uniqueness of $u$-measures and SRB measures}\label{srb:uniq} 
In the case of a hyperbolic attractor, topological transitivity of $f|\Lambda$ guarantees that there is a unique $u$-measure for $f$ on $\Lambda$. In contrast, in the partially hyperbolic situation, even topological mixing is not enough to guarantee that there is a unique $u$-measure. Indeed, consider $F=f_1\times f_2$, where $f_1$ is a topologically transitive Anosov diffeomorphism and $f_2$ a diffeomorphism close to the identity. Then 
$F$ is partially hyperbolic, and any measure $\mu=\mu_1\times \mu_2$, where $\mu_1$ is the unique SRB measure for $f_1$ and $\mu_2$ any $f_2$-invariant measure, is a $u$-measure for $F$. Thus, $F$ has a unique $u$-measure if and only if $f_2$ is uniquely ergodic. On the other hand, $F$ is topologically mixing if and only if $f_2$ is topologically mixing. 
\begin{theorem}[\cite{BDPP}]\label{bdpp1}
Let $f$ be a $C^{1+\alpha}$ diffeomorphism of a compact smooth Riemannian manifold $M$ possessing a partially hyperbolic attractor $\Lambda$. Assume that: 
\begin{enumerate}
\item there exists a $u$-measure $\mu$ for $f$ with respect to which $f$ has negative central exponents on an invariant subset $A\subset\Lambda$ of positive $\mu$-measure; 
\item for every $x\in\Lambda$ the orbit of the global strongly unstable manifold $W^u (x)$ is dense in $\Lambda$. 
\end{enumerate}
Then $\mu$ is the only $u$-measure for $f$ and $f$ has negative central  exponents at 
$\mu$-almost every $x\in\Lambda$. In particular, $(f,\mu)$ is ergodic, 
$\text{supp}(\mu)=\Lambda$, and the basin $B_\mu$ has full volume in the topological basin of attraction of $\Lambda$. $\mu$ is the only SRB measure for $f$.
\end{theorem}

Let us comment on the assumption of this theorem. Shub and Wilkinson \cite{ShWi} considered the direct product $F_0=f\times\text{Id}$, where $f$ is a linear Anosov diffeomorphism and the identity acts on the circle. The map $F_0$ preserves volume. They showed that arbitrary close to $F_0$ (in the $C^1$ topology) there is a volume-preserving diffeomorphism $F$ whose only central exponent is negative on the whole of $M$. The result continues to hold for any small perturbation of $F$.

Bonatti and Diaz \cite{BoDi} have shown that there is an open set of transitive diffeomorphisms near $F_0=f\times\text{Id}$ ($f$ is an Anosov diffeomorphism and $\text{Id}$ is the identity map of any manifold) as well as near the time-$1$ map $F_0$ of a topologically transitive Anosov flow. This result was used by Bonatti, Diaz and Ures \cite{BDU} to construct examples of partially hyperbolic systems with minimal unstable foliation (i.e., every unstable leaf is dense in the manifold itself). 

If $f$ is a small perturbation of $F_0$ then $f$ is partially hyperbolic and by \cite{HPS2}, the central distribution of $f$ is integrable. Furthermore, the central leaves are compact in the first case (when $F_0=f\times\text{Id}$) and there are compact leaves in the second case (when $F_0$ is the time-$1$ map of a topologically transitive Anosov flow). It is shown in 
\cite{BDPP} that if there is a compact periodic central leaf $C$ for $f$ such that $f^n (C)=C$ and the restriction $f^n|C$ is a minimal transformation, then the unstable foliation for $f$ is minimal. 

\subsection{Stable ergodicity for dissipative maps} Let $\Lambda_f$ be a topological attractor for a diffeomorphism $f$. We say that $f$ is \emph{stably ergodic} if there exists a neighborhood $\mathcal U$ of $f$ in $\text{Diff}^r(M)$, $r\ge 1$ such that  any diffeomorphism $g\in{\mathcal U}$ possesses a topological attractor $\Lambda_g$ and there is a unique SRB measure $\mu_g$ on $\Lambda_g$ (and hence, $g$ is ergodic with respect to $\mu_g$). This is an analog of the notion of stable ergodicity of systems preserving a given smooth measure, which was introduced by Pugh and Shub, \cite{PuSh}. For systems with topological attractors smooth measures are replaced by (unique) SRB measures.

If the attractor $\Lambda_f$ is (partially) hyperbolic then there exists a neighborhood 
$\mathcal U$ of $f$ in $\text{Diff}^1(M)$ such that any diffeomorphism $g\in{\mathcal U}$ possesses a (partially) hyperbolic attractor $\Lambda_g$. 

\begin{theorem}[\cite{BDPP}]\label{bdpp2}
Let $\Lambda_f$ be a partially hyperbolic attractor for a diffeomorphism $f$. If $f$ satisfies the conditions of Theorem \ref{bdpp1}, then $f$ is stably ergodic with $r=1+\alpha$. 
\end{theorem}

The stable ergodicity of partially hyperbolic attractors with positive central exponents was studied by V\'asquez \cite{Vas} who proved a result similar to Theorem~\ref{bdpp2} under the stronger requirement that there is a {\bf unique} $u$-measure with positive central exponents on a subset of {\bf full} measure.

\section{SRB measures for attractors with dominated splitting}\label{dom-split}

Let $f$ be a $C^{1+\alpha}$ diffeomorphism and $A$ a forward-invariant compact set. 
A splitting 
\[
T_AM=E^{cs}\oplus E^{cu}
\] 
is \emph{dominated} if there is $\chi<1$ such that
$$
\|df|_{E^s(x)}\|<\chi\|df|_{E^u(x)}^{-1}\|^{-1}
$$
for all \( x\in A \). The set 
$$
\Lambda=\bigcap_{j=0}^\infty f^j(A)
$$
is an \emph{attractor for $f$ with dominated splitting}.

\begin{remark}
We note that the domination condition does not imply anything about whether the derivative is contracting or expanding in each of the sub-bundles, but just that it is contracting in \( E^{cs} \) \emph{relative} to \( E^{cu} \).  There are no general results on the existence of SRB measures for dominated splitting in full generality, but only with some additional assumptions, including that either \( E^{cs} \) is uniformly contracting or that \( E^{cu} \) is uniformly expanding. It will be useful to distinguish these two cases by using the notation 
\begin{equation}\label{domsplit}
T_AM=E^{s}\oplus E^{cu} \quad \text{ and } \quad T_AM=E^{cs}\oplus E^{u}
\end{equation}
respectively. Notice that if \( E^{cs} \) is uniformly contracting and \( E^{cu} \) is uniformly expanding then we are in the uniformly hyperbolic situation with a splitting  \( T_AM=E^{s}\oplus E^{u} \). A partially hyperbolic splitting \( T_{A}M=E^{s}\oplus E^{c}\oplus E^{u} \) is also a special case of dominated decomposition since we can combine the central bundle with either the stable one to get a dominated splitting of the form \( T_{A}M=E^{cs}\oplus E^{u} \) where \( E^{cs}=E^{s}\oplus E^{c} \) or combine the central bundle with the unstable one to get a splitting of the form \( T_{A}M=E^{s}\oplus E^{cu} \) where \( E^{cu}=E^{c}\oplus E^{u} \). We emphasize however that dominated decompositions, even under the additional assumptions that one of the two bundles is uniform, does not imply that the center-stable bundle \( E^{cs} \) or the center-unstable bundle \( E^{cu} \) can be further split into a center bundle \( E^{c} \) and either a stable or unstable bundle. This additional splitting can play a key role in some results about partially hyperbolic attractors which therefore do not immediately extend to attractors with dominated splittings. 
\end{remark}

The existence of SRB measures has been proved for attractors with dominated splittings of the form \eqref{domsplit} as long there is also some degree of (non-uniform) contraction or expansion in the center-stable or center-unstable bundles respectively.  

\begin{theorem}[\cite{BV}]\label{th:BV}
Let \( f\colon M\to M \) be a  \( C^{1+\alpha} \) diffeomorphism, $A$ a forward invariant compact set with \( m(A)>0 \) on which \(  f  \) admits a dominated splitting of the form 
\[
T_{A}M=E^{cs}\oplus E^{u}
\]
and suppose that the dynamics on \( E^{cs} \) is \emph{nonuniformly contracting}:\begin{equation} \label{NUE1}
\limsup_{n\to+\infty} \frac{1}{n}
\log \|df^{n}|{E^{cs}_{x}}\| <0
\end{equation}
for all \(  x\in A  \). Then \( f \) has finitely many ergodic SRB measures and the union of their basins has full Lebesgue measure in the topological basin of \( A \).
\end{theorem}

This theorem can be proved by a ``push-forward'' argument that essentially follows the lines of the proof of Theorems~\ref{existence-srb} and \ref{bdpp1} in the uniformly (completely or partially) hyperbolic case. The full definition of SRB measure requires some (possibly nonuniform) hyperbolicity (i.e., contraction and expansion in all directions) in order to guarantee existence and absolute continuity of the stable foliation which in turns guarantees that the measure is a physical measure. Condition \eqref{NUE1} bridges this gap and yields the statement of Theorem~\ref{th:BV}. 

\begin{theorem}[\cite{ABV, ADLP}]\label{th:ADLP}
Let \( f\colon M\to M \) be a  \( C^{1+\alpha} \)
diffeomorphism,  $A$   a forward invariant compact set with \( m(A)>0 \) on which
\(  f   \) admits a dominated splitting of the form  
\[
T_{A}M=E^{s}\oplus E^{cu}
\]
and suppose that the dynamics on  \( E^{cu} \) is \emph{nonuniformly expanding}:
\begin{equation} \label{NUE2}
\limsup_{n\to+\infty} \frac{1}{n}
\sum_{j=1}^{n} \log 
(\|df^{-1}|E^{cu}_{f^j(x)}\|^{-1})
 >  \eps
\end{equation}
for some \( \eps>0 \) and for all \(  x\in A  \). 
Then \( f \) has finitely many ergodic SRB measures and the union of their basins has full Lebesgue measure in the topological basin of \( A \). 
\end{theorem}
This result was first proved in \cite{ABV} under the slightly stronger assumption obtained by replacing \( \limsup \) by \( \liminf \) in \eqref{NUE2}. Because the expansion along $E^{cu}$ is nonuniform, the argument there uses a more sophisticated version of the geometric ``push-forward'' argument of Lebesgue measure \( m_{V^{u}(x)} \) on the local unstable manifolds than the one outlined in Section \ref{exist-hyper-attr1}. 

Indeed, it is not longer true that for every \( k \) we can divide \( V^{u}(x) \) into pieces each of which grows to large scale with bounded distortion at time \( k \). Instead this will be true just for \emph{some} points in \(V^{u}(x) \), precisely those for which \( k \) is a \emph{hyperbolic time}.  The images at time \( k \) of other parts of \( V^{u}(x) \) may be very small and/or very distorted. In particular it is no longer the case that \( f^{k}_{*}(V^{u}(x)) \) is supported on a collection of uniformly large unstable disks. Nevertheless some points do eventually have hyperbolic times and therefore it is the case that  \emph{some part} of the measures \( f^{k}_{*}(V^{u}(x)) \), and therefore some part of the measures \( \nu_{n} \), are supported on some such collection of uniformly large unstable disk. Thus it is possible to write  the measures \( \nu_{n} \) as 
 \begin{equation*}
\nu_{n}=\nu_{n}'+\nu_{n}''
\end{equation*}
where \( \nu_{n}' \) is the ``good'' part of the measure supported on a collection of uniformly large unstable disks and \( \mu_{n}'' \) is the ``bad'' part on which we have little control. The strengthened version of condition \eqref{NUE2}, with a \( \liminf \) instead of a \( \limsup \),  implies that almost every point has a positive density of hyperbolic times and thus  the good part of the measure \( \mu_{n}' \) forms a proportion of the overall measure \( \mu_{n} \) that is uniformly bounded below in \( n \) and  it is therefore possible to essentially recover a version of the original argument of Sinai, Ruelle, Bowen and show that there exists a limit measure \( \mu' \) which has absolutely continuous conditional measures and is therefore an SRB measure (see more details in the outline of the proof of Theorem \ref{thm:7.1}).  

Replacing the \emph{liminf} by a \emph{limsup} as stated in condition \eqref{NUE2} still implies an \emph{infinite number} of hyperbolic times and this still allows us to split the measures into \( \nu'_{n} \) and \( \nu''_{n} \), but does not imply a \emph{positive density} of hyperbolic times  and thus makes it impossible to obtain a uniform lower bound for the mass of the measures \( \nu'_{n} \) and to complete the proof using the natural ``push-forward'' argument. The full proof of Theorem \ref{th:ADLP} is thus obtained in \cite{ADLP} using the \emph{inducing} or \emph{Young tower} approach mentioned above. 
 In certain respects there are of course still some similarities with the classical approach in the sense that the Young tower structure also relies on constructing some region where large unstable disks accumulate, and defining an induced map on this region. The problem of the low asymptotic frequency of hyperbolic times does not disappear in this approach but is rather ``translated'' into the problem of integrability of the return times, which can be resolved by using a different kind of argument. 
 
\section{SRB measures for non-uniformly hyperbolic attractors}\label{non-uni-attr}

If $\mu$ is an SRB measure, then every point in the positive Lebesgue measure set 
$B_\mu$ has non-zero Lyapunov exponents. A natural and interesting question is whether the converse holds true, in essence formulated in the following conjecture by Viana, \cite{Viana}:
\begin{conjecture}
If a smooth map has only non-zero Lyapunov exponents at Lebesgue almost every point, then it admits an SRB measure. 
\end{conjecture}
The results of the previous section can be viewed as partial progress in the direction of Viana's conjecture by proving the existence of SRB measures under the assumptions of non-zero Lyapunov exponents and additional conditions that the system has a dominated splitting and either stable or unstable direction is uniformly hyperbolic.\footnote{Although note that the `mostly expanding' condition of Theorem \ref{th:ADLP} is slightly more restrictive than saying that all Lyapunov exponents in $E^{cu}$ are positive.}  The presence of the dominated splitting means that one does not need to worry too much about the geometry of the stable and unstable manifolds, and only needs to take care of the expansion and contraction properties.

In general, however,  the geometric properties of the system are not uniform. In a ``fully'' non-uniformly hyperbolic system the splitting $E^s\oplus E^u$ is only measurable, and the angle between the stable and unstable subbundles is arbitrarily small. In this section we describe two general results that apply in this setting.

The first significant results on SRB measures for non-uniformly hyperbolic systems were those for the attractors for certain special parameters of the H\'enon family of maps obtained in \cite{BenCar}. These attractors have a fully nonuniformly hyperbolic structure which can be described relatively explicitly and, taking advantage of several specific characteristics of this structure, an SRB measure for the attractors was constructed first in \cite{BY1} using the push-forward argument in what is essentially a further variation of argument using the splitting \( \nu_{n}=\nu'_{n}+\nu''_{n} \) described above in the case of dominated splittings, and later in \cite{BY2} using the induced map argument in a construction which effectively became the model case study for Young's full development of this general approach in~\cite{You98}. 

We describe here two recent results which develop a general framework for the construction of SRB measures for fully general nonuniformly hyperbolic attractors. The first one applies in arbitrary dimension and gives the existence of SRB measures under some slight strengthening of usual  nonuniform hyperbolicity. The construction in this case is, once again, a further and even more sophisticated refinement of the ``push-forward'' argument described above. The second result applies only in dimension two but proves the existence of SRB measure under no more than standard non-uniform hyperbolicity conditions plus a natural recurrence assumption. The proof of this result is instead based on the Young tower approach and indeed a corollary of independent interest, is that in dimension two, \emph{any} SRB measure can in principle be obtained through a Young tower.  

\subsection{Effective hyperbolicity} We make the following standing assumption.
\begin{enumerate}[label=\textup{\textbf{(H)}}]
\item\label{H} There exists a forward-invariant set $A\subset U$ of positive volume with two measurable cone families $K^s(x), K^u(x)\subset T_xM$ such that
\begin{enumerate}
\item $\overline{Df(K^u(x))}\subset K^u(f(x))$ for all $x\in A$;
\item $\overline{Df^{-1}(K^s(f(x)))}\subset K^s(x)$ for all $x\in f(A)$.
\item $K^s(x)=K(x,E^s(x),a_s(x))$ and $K^u(x)=K(x,E^u(x),a_u(x))$ are such that 
$T_xM=E^s(x)\oplus E^u(x)$; moreover $d_s=\dim E^s(x)$ and $d_u=\dim E^u(x)$ do not depend on $x$.
\end{enumerate}
\end{enumerate}
Such cone families automatically exist if $f$ is uniformly hyperbolic on $\Lambda$. We emphasize, however, that in our setting $K^{s,u}$ are not assumed to be continuous, but only measurable and the families of subspaces $E^{u,s}(x)$ are not assumed to be invariant. 

Let $A\subset U$ be a forward-invariant set satisfying \ref{H}. Define 
$$
\begin{aligned}
\lambda^u(x) &= \inf \{\log \|Df(v)\| \mid v\in K^u(x), \|v\|=1 \}, \\
\lambda^s(x) &= \sup \{\log \|Df(v)\| \mid v\in K^s(x), \|v\|=1 \}.
\end{aligned}
$$
Note that if the splitting $E^s\oplus E^u$ is dominated, then we have 
$\lambda^s(x)<\lambda^u(x)$ for every $x$. Thus we define the \emph{defect from domination} at $x$ to be
$$
\Delta(x) = \tfrac 1\alpha \max(0,\ \lambda^s(x) - \lambda^u(x)),
$$
where $\alpha\in (0,1]$ is the H\"older exponent of $Df$. Roughly speaking, $\Delta(x)$ controls how much the curvature of unstable manifolds can grow as we go from $x$ to 
$f(x)$.

The following quantity is positive whenever $f$ expands vectors in $K^u(x)$ and contracts vectors in $K^s(x)$:  
$$
\lambda(x) = \min ( \lambda^u(x) - \Delta(x), -\lambda^s(x)).
$$
The \emph{upper asymptotic density} of $\Gamma\subset\mathbb{N}$ is
$$
\overline{\delta}(\Gamma)=\limsup_{N\to\infty}\frac1N\#\big(\Gamma\cap [0,N)\big).
$$
An analogous definition gives the lower asymptotic density $\underline{\delta}(\Gamma)$.

Denote the angle between the boundaries of $K^s(x)$ and $K^u(x)$ by
$$
\theta(x)=\inf\{\measuredangle(v,w)\colon v\in K^u(x), w\in K^s(x) \}.
$$
We say that a point $x\in A$ is \emph{effectively hyperbolic} if 
\begin{gather}
\tag{\textbf{EH1}}\label{EH1}
\liminf_{n\to\infty}\frac1n\sum_{k=0}^{n-1}\lambda(f^k(x))>0,\\
\tag{\textbf{EH2}}\label{EH2}
\lim_{\bar\theta\to 0}\overline\delta\{n\mid\theta(f^n(x))<\bar\theta\}=0.
\end{gather}
Condition \eqref{EH1} says that not only are the Lyapunov exponents of $x$ positive for vectors in $K^u$ and negative for vectors in $K^s$, but $\lambda^u$ gives enough expansion to overcome the `defect from domination' given by $\Delta$.

Condition \eqref{EH2} requires that  the frequency with which the angle between the stable and unstable cones drops below a specified threshold $\bar\theta$ can be made arbitrarily small by taking the threshold to be small.  

If $\Lambda$ is a hyperbolic attractor for $f$, then {\bf every} point $x\in U$ is effectively hyperbolic, since there are $\bar\lambda,\bar\theta>0$ such that 
$\lambda^s(x)\le -\bar\lambda$, $\lambda^u(x)\ge\bar\lambda$, and 
$\theta(x)\ge\bar\theta$ for every $x\in U$, so that $\Delta(x)=0$ and 
$\lambda(x)\ge\bar\lambda$.

Let $A$ satisfy \ref{H}, and let $S\subset A$ be the set of effectively hyperbolic points.
Observe that effective hyperbolicity is determined in terms of a forward asymptotic property of the orbit of $x$, and hence $S$ is forward-invariant under $f$. The following result is proved in \cite{CDP}.

\begin{theorem}\label{thm:7.1}
Let $f$ be a $C^{1+\alpha}$ diffeomorphism of a compact manifold $M$, and $\Lambda$ a topological attractor for $f$.  Assume that 
\begin{enumerate}
\item $f$ admits measurable invariant cone families as in \ref{H}; 
\item the set $S$ of effectively hyperbolic points satisfies $m(S)>0$.
\end{enumerate}
Then $f$ has an SRB measure supported on $\Lambda$.
\end{theorem}
A similar result can be formulated given information about the set of effectively hyperbolic points on a single `approximately unstable' submanifold usually called \emph{admissible}. The set of admissible manifolds that we will work with is related to $\mathcal{R}_{\bf I}$ from Section \ref{hyp-attr}, but the precise definition is not needed for the statement of the theorem; all we need here is to have $T_x W\subset K^u(x)$ for `enough' points $x$. 
$W\subset U$. Let $d_u$, $d_s$, and $A$ be as in \ref{H}, \eqref{EH1} and \eqref{EH2}, and let $W\subset U$ be an embedded submanifold of dimension $d_u$.

\begin{theorem}\label{thm:CDP}
Let $f$ be a $C^{1+\alpha}$ diffeomorphism of a compact manifold $M$, and $\Lambda$ a topological attractor for $f$.  Assume that 
\begin{enumerate}
\item $f$ admits measurable invariant cone families as in \ref{H};
\item there is a $d_u$-dimensional embedded submanifold $W\subset U$ such that 
$m_W(\{x\in S\cap W\mid T_xW\subset K^u(x)\})>0$.
\end{enumerate}
Then $f$ has an SRB measure supported on $\Lambda$.
\end{theorem}

We outline the proof of this statement to illustrate the geometric approach in the settings of non-uniformly hyperbolic attractors. We follow the same ideas as in Section \ref{hyp-attr}, but there are two major obstacles to overcome.
\begin{enumerate}
\item The action of $f$ along admissible manifolds is not necessarily uniformly expanding.
\item Given $n\in\NN$ it is no longer necessarily the case that $f^n(W)$ contains any admissible manifolds in $\RRR_{\bf I}$, let alone that it can be covered by them.  When $f^n(W)$ contains some admissible manifolds, we will need to control how much of  it  can be covered.
\end{enumerate}
To address the first of these obstacles, we need to consider admissible manifolds for which we control not only the geometry but also the dynamics; thus we will replace the collection 
$\RRR_\KK'$ from before with a more carefully defined set (in particular, $\KK$ will include more parameters). Since we do not have uniformly transverse invariant subspaces 
$E^{u,s}$, our definition of an admissible manifold also needs to specify which subspaces are used, and the geometric control requires an assumption about the angle between them.  

Given $\theta,\gamma,\kappa,r>0$, write $\II=(\theta,\gamma,\kappa,r)$ and consider the following set of ``$(\gamma,\kappa)$-admissible manifolds of size $r$ with transversals controlled by $\theta$'':
\begin{multline}\label{eqn:PPP}
\PPP_\II=\{\exp_x(\graph\psi)\mid x\in\overline{f(U)},\ T_x M=G\oplus F,\
G \subset \overline{K^u(x)},\\ 
\angle(G,F)\ge\theta,\ \psi\in C^{1+\alpha}(B_G(r), F)\text{ satisfies~\eqref{admissibe-man}}\}.
\end{multline}
Elements of $\PPP_\II$ are admissible manifolds with controlled geometry.  We also impose a condition on the dynamics of these manifolds.  Fixing $C, \bl>0$, write 
$\JJ=(C,\bl)$ and consider for each $N\in\NN$ the collection of sets
\begin{multline}\label{eqn:QQQ}
\QQQ_{\JJ,N}=\{f^N(V_0)\mid V_0\subset U, \text{ and for every } y,z\in V_0, \text{ we have} \\
d(f^j(y), f^j(z))\le Ce^{-\bl(N- j)}d(f^N(y),f^N(z))\text{ for all } 0\le j\le N\}.
\end{multline} 
Elements of $\PPP_\II\cap \QQQ_{\JJ,N}$ are admissible manifolds with controlled geometry and dynamics in the unstable direction. We also need a parameter $\beta>0$ that controls the dynamics in the stable direction, and another parameter $L>0$ that controls densities in standard pairs. Then writing $\KK=\II\cup\JJ\cup\{\beta,L\}$, we obtain a set $\RRR_{\KK,N}\subset\PPP_\II\cap\QQQ_{\JJ,N}$ for which we have the added restriction that we control the dynamics in the stable direction; the corresponding set of standard pairs is written $\RRR_{\KK,N}'$.  

The set $\RRR'_{\KK,N}$ carries a natural product topology in which $\RRR_{\KK,N}'$ is compact and the map $\Phi$ defined in \eqref{eqn:Phi2} is continuous.

As before, let $\MMM_{\le 1}(\RRR'_{\KK,N})$ denote the space of measures on 
$\RRR'_{\KK,N}$ with total weight at most 1. The resulting space of measures on $U$ plays a central role:
\begin{equation}\label{eqn:Mach}
\MMM_{\KK,N}=\Phi(\MMM_{\leq 1}(\RRR'_{\KK,N})).
\end{equation}
Measures in $\MMM_{\KK,N}$ have uniformly controlled geometry, dynamics, and densities via the parameters in $\KK$, and $\MMM_{\KK,N}$ is compact. However, at this point we encounter the second obstacle mentioned above: because $f(W)$ may not be covered by admissible manifolds in $\RRR_{\KK,N}$, the set $\MMM_{\KK,N}$ is not 
$f_*$-invariant.  

To address this, one must establish good recurrence properties to $\MMM_{\KK,N}$ under the action of $f_*$ on $\MMM(\overline{U})$; this can be done via effective hyperbolicity.

Consider for $x\in A$ and $\bl>0$ the set of \emph{effective hyperbolic times}
\begin{equation}\label{eqn:ehtimes}
\Gamma^e_\bl(x)=\bigg\{n\mid\sum_{j=k}^{n-1}(\lambda^u-\Delta)(f^j(x))\ge\bl(n-k) \text{ for all }0\le k<n\bigg\}.
\end{equation}
Results from \cite{CP} show that for every $x$ and almost every effective hyperbolic time $n\in \Gamma^e_\bl(x)$, there is a neighborhood 
$W_n^x\subset W$ containing $x$ such that $f^n(W_n^x)\in\PPP_\II\cap\QQQ_{\JJ,N}$.  With a little more work, one can produce a ``uniformly large'' set of points $x$ and times $n$ such that $f^n(W_n^x)\in\RRR_{\KK,N}$, and in fact $f_*^n m_{W_n^x}\in\MMM_{\KK,N}$. 
Then this can be used to obtain measures $\nu_n \in\MMM_{\KK,N}$ such that 
\begin{equation}\label{eqn:nun-ulim}
\nu_n\le\mu_n = \tfrac 1n \textstyle\sum_{k=0}^{n-1} f_*^k m_W
\qquad\text{and}\qquad
\ulim_{n\to\infty} \|\nu_n\| > 0.
\end{equation}
Once this is achieved, compactness of $\MMM_{\KK,N}$ guarantees existence of a non-trivial $\nu\in\bigcap_N\MMM_{\KK,N}$ such that $\nu\le\mu=\lim_k \mu_{n_k}$. In order to apply the absolute continuity properties of $\nu$ to the measure $\mu$, one must define a collection $\Mac$ of measures with good absolute continuity properties along admissible manifolds, for which there is a version of the Lebesgue decomposition theorem that gives 
$\mu = \mu^{(1)} + \mu^{(2)}$, where $\mu^{(1)}\in\Mac$ is invariant. This measure is non-trivial since $0\ne\nu\le\mu^{(1)}$, and the definition of $\RRR'_{\KK,N}$ guarantees that the set of points with non-zero Lyapunov exponents has positive measure with respect to 
$\nu$, and hence also with respect to $\mu^{(1)}$. Thus some ergodic component of 
$\mu^{(1)}$ is hyperbolic, and hence is an SRB measure.  

\subsection{Maps on the boundary of Axiom A: neutral fixed points}

We give a specific example of a map for which the conditions of Theorem \ref{thm:CDP} can be verified. Let $f\colon U\to M$ be a $C^{1+\alpha}$ Axiom A diffeomorphism onto its image with $\overline{f(U)}\subset U$, where $\alpha\in (0,1)$. Suppose that $f$ has one-dimensional unstable bundle.

Let $p$ be a fixed point for $f$. We perturb $f$ to obtain a new map $g$ that has an indifferent fixed point at $p$. The case when $M$ is two-dimensional and $f$ is volume-preserving was studied by Katok. We allow manifolds of arbitrary dimensions and (potentially) dissipative maps. For example, one can choose $f$ to be the Smale--Williams solenoid or its sufficiently small perturbation.  

We describe a specific perturbation of $f$ for which the conditions of the main theorem can be verified; one can also describe a general set of conditions on the return map through the region of the perturbation \cite[Theorem 2.3]{CDP}. We suppose that there exists a neighborhood $Z\ni p$ with local coordinates in which $f$ is the time-$1$ map of the flow generated by
$$
\dot x = Ax
$$
for some $A\in GL(d,\mathbb{R})$. Assume that the local coordinates identify the splitting $E^u\oplus E^s$ with $\mathbb{R}\oplus\mathbb{R}^{d-1}$, so that $A=A_u\oplus A_s$, where $A_u=\gamma\text{Id}_u$ and $A_s=-\beta\text{Id}_s$ for some $\gamma,\beta>0$.  In the Katok example we have $d=2$ and $\gamma=\beta$ since the map is 
area-preserving. 

Now we use local coordinates on $Z$ and identify $p$ with $0$. Fix $0<r_0<r_1$ such that $B(0,r_1)\subset Z$, and let $\psi\colon Z\to [0,1]$ be a $C^{1+\alpha}$ function such that
\begin{enumerate}
\item $\psi(x)=\|x\|^\alpha$ for $\|x\|\le r_0$;
\item $\psi(x)=1$ for $\|x\|\ge r_1$;
\item $\psi(x)>0$ for $x\ne 0$ and $\psi'(x)>0$.
\end{enumerate}
Let $\mathcal{X}\colon Z\to\mathbb{R}^d$ be the vector field given by 
$\mathcal{X}(x)=\psi(x) Ax$. Let $g\colon U\to M$ be given by the time-$1$ map of this vector field on $Z$ and by $f$ on $U\setminus Z$. Note that $g$ is $C^{1+\alpha}$ because $\mathcal{X}$ is $C^{1+\alpha}$. In \cite{CDP} it is shown that $g$ satisfies the conditions of Theorem \ref{thm:CDP} (in fact, a slightly more general version of this theorem), which proves the following.
\begin{theorem}\label{thm:katok}
The map $g$ has an SRB measure.
\end{theorem}
Note that $g$ does not have a dominated splitting because of the indifferent fixed point.
We also observe that if $\psi$ is taken to be $C^\infty$ away from $0$, then $g$ is also $C^\infty$ away from the point $p$.  

\subsection{Two-dimensional non-uniformly hyperbolic attractors}

Let $f$ be a $C^{1+\alpha}$ diffeomorphism of a compact surface which is non-uniformly hyperbolic on an invariant set $\Lambda$.\footnote{We stress that non-uniform hyperbolicity on an invariant set does not require presence of any invariant measure. On the other hand if $f$ preserves a hyperbolic measure then there is an invariant set $\Lambda$ on which $f$ is non-uniformly hyperbolic.}  Let $\Lambda_\ell$ be a regular set.
%

\begin{definition}\label{def:rectangles}
A subset  \(  \mathcal R \subset \Lambda \)  is a  \emph{rectangle} if \( \mathcal R \subseteq \Lambda_{\ell} \) for some \( \ell \) and   for every $x,y\in\mathcal R$,   $[x,y] := V^s(x)\cap V^u(y)$ consists of  a single point and \( [x,y]\in \mathcal R \).
\end{definition}
In other words the set $\mathcal R$ has \emph{hyperbolic product structure}. 
Notice that  a single point \( x\in \Lambda \) can be considered a (trivial) rectangle, but we in general we will always assume that our rectangle consists of an uncountable set of points, in which case \( \mathcal R \) contains lots of subrectangles. Indeed for any two distinct points \( p, q\in \mathcal R \),  their stable and unstable curves bound an open domain \( \widehat {\mathcal R}_{p,q}\subset  M \) and the intersection of \( \mathcal R \) with the closure of \( \widehat{\mathcal R}_{p,q} \) is also a rectangle, which we denote by \( \mathcal R_{p,q} \). 

\begin{definition}\label{def:rectprop}
A rectangle \( \mathcal R \) is 
\begin{enumerate}
\item \emph{nice} if \( \mathcal R = \mathcal R_{p,q} \) for some periodic points \( p,q \);
\item
\emph{fat} if there exists some \( x\in \mathcal R \) with \( m_{V^u(x)}(\mathcal R\cap V^{u}(x))>0 \); 
\item 
\emph{recurrent} if  every \( x\in \mathcal R \) returns to \( \mathcal R \)  in the future and in the past. 
\end{enumerate}
\end{definition}
We remark that the existence of rectangles with the \emph{nice} property is easily satisfied in two dimensions for any hyperbolic set \( \Lambda \) as above. The other two conditions are elementary and natural conditions which are easily seen to be necessary for the existence of an SRB measure. It is proved in \cite{CLP} that in fact these are also sufficient conditions. 

\begin{theorem}[\cite{CLP}]\label{thm:nonunifSRB}
The map \( f \) has an SRB measure if and only if it admits a nice fat recurrent rectangle. 
\end{theorem}
As mentioned above, the proof of this result relies on the construction of a Young tower, which is in itself a result of independent interest. Indeed, the technical core of the proof of Theorem \ref{thm:nonunifSRB} is a result which gives some extremely simple and natural conditions for the existence of a Young tower in the two dimensional case. To formulate this result we give here the definition of the topological structure of a Young tower, which is the most difficult part of the construction and refer the reader to \cite{CLP} or the original paper of Young \cite{You98} for the full definition which includes some hyperbolicity and distortion conditions.

\begin{definition} Let \( \Gamma \) be a rectangle. 
\begin{enumerate}
\item $\Gamma^s \subset\Gamma$ is an \(s\)\emph{-subset} of \( \Gamma \) if $x\in\Gamma^s$ implies $V^s(x)\cap\Gamma\subset\Gamma^s$;
\item $\Gamma^u \subset\Gamma$ is a \(u\)\emph{-subset} of \( \Gamma \) if $x\in\Gamma^u$ implies $V^u(x)\cap\Gamma\subset\Gamma^u$.
\end{enumerate}
\end{definition}  
  
Let \( \tau\colon \Gamma\to \mathbb N \) be the first return time function of points of \( \Gamma \) to \( \Gamma \) (which is defined at all points of \( \Gamma \) if \( \Gamma \) is recurrent) and let \(  \mathbb N_\tau \) denote the set of  times which occur as first return times of points in  \( \Gamma \).   Let \( \Gamma \) be a recurrent rectangle. 

\begin{definition}\label{def:FRYT}
\( \Gamma \) supports a \emph{Topological First Return Young Tower}~if
$$
\Gamma^S_i :=\{x\in\Gamma: \tau(x)=i\} \quad \text{ and  } \quad  \Gamma^U_i:=f^i(\Gamma^S_i)
$$
are \( s\)-subsets and \( u\)-subsets respectively of $\Gamma$ for every $i\in\mathbb N_\tau$.
\end{definition}

\begin{theorem}[\cite{CLP}]\label{thm:maintop}
Let \( \Gamma_{0} \)  be a nice recurrent rectangle.  Then there exists a nice recurrent rectangle \( \Gamma\supset \Gamma_{0} \) which supports  a Topological First Return Young Tower. 
\end{theorem}

The property that \( \Gamma\supset \Gamma_{0} \) implies that if \( \Gamma_{0} \) is a fat rectangle then the same holds for \( \Gamma \) and thus an immediate consequence of Theorem \ref{thm:maintop} is that if \( \Gamma_{0} \) is a nice fat recurrent rectangle then there exists another, also nice fat recurrent, rectangle which supports a First return Topological Young Tower. It is then possible to prove that the required hyperbolicity and distortion estimates are satisfied  which imply the existence of an SRB measure by \cite{You98}, thus obtaining Theorem \ref{thm:nonunifSRB}. 

\section{SRB measures for hyperbolic attractors with singularities}
\label{att-sing}

\subsection{Topological attractors with singularities} 
Let $M$ be a smooth compact manifold, $U\subset M$ an open bounded connected subset, \emph{the trapping region}, $N\subset U$ a closed subset and 
$f\colon U\setminus N\to U$ a $C^2$ diffeomorphism such that 
\begin{equation}\label{singularities}
\begin{aligned}
\|d^2f_x\|&\le C_1d(x,\mathcal{S}^+)^{-\alpha_1} \text{ for any } x\in U\setminus N,\\
\|d^2f^{-1}_x\|&\le C_2d(x,\mathcal{S}^-)^{-\alpha_2} \text{ for any } x\in f(U\setminus N),
\end{aligned}
\end{equation}
where $\mathcal{S}^+=N\cup\partial U$ is the \emph{singularity set for $f$} and 
$\mathcal{S}^-=f(\mathcal{S}^+)$ that is
$$
\mathcal{S}^-=\{y\in U:\text{ there is }z\in\mathcal{S}^+\text{ and }z_n\in U\setminus \mathcal{S}^+ \text{ such that } z_n\to z, f(z_n)\to f(z)\}
$$
is the \emph{singularity set for $f^{-1}$}. We will assume that 
$m(\mathcal{S}^+)=m(\mathcal{S}^-)=0$.

Define
$$
U^+=\{x\in U: f^n(x)\notin \mathcal{S}^+, n=1,2,\dots\}
$$
and the \emph{topological attractor with singularities}  
$$
D=\bigcap_{n\ge 0}f^n(U^+), \quad \Lambda=\bar D.
$$
Given $\varepsilon>0$ and $\ell>1$, set
$$
\begin{aligned}
D^+_{\varepsilon,\ell}&=\{z\in\Lambda : d(f^n(z),\mathcal{S}^+)\ge \ell^{-1}e^{-\varepsilon n}, n=0,1,2,\dots\},\\
D^-_{\varepsilon,\ell}&=\{z\in\Lambda : d(f^n(z),N^-)\ge \ell^{-1}e^{-\varepsilon n}, n=0,1,2,\dots\},\\
D^0_{\varepsilon,\ell}&=D^+_{\varepsilon,\ell}\bigcap D^-_{\varepsilon,\ell},\\
D^0_\varepsilon&=\bigcup_{\ell\ge 1} D^0_{\varepsilon,\ell}.
\end{aligned}
$$
The set $D^0_\varepsilon$ is the \emph{core} of the attractor and it may be an empty set as it may be the set $D$.

\begin{theorem}[\cite{Pes2}]
Assume that there are $C>0$ and $q>0$ such that for any $\varepsilon>0$ and $n>0$
\begin{equation}\label{core}
m(f^{-n}(\mathcal{U}(\varepsilon,\mathcal{S}^+)\cap f^n(U^+)))\le C\varepsilon^q,
\end{equation}
where $\mathcal{U}(\varepsilon,\mathcal{S}^+)$ is a neighborhood of the (closed) set $\mathcal{S}^+$. Then there is an invariant measure $\mu$ on $\Lambda$ such that, $\mu(D^0_\varepsilon)>0$, in particular, the core is not empty.
\end{theorem}

\subsection{Hyperbolic attractors with singularities} We say that a topological attractor with singularities $\Lambda$ is \emph{hyperbolic}, if there exist two families of stable and unstable cones 
$$
K^s(x)=K(x,E_1(x),\theta(x)), \, K^u(x)=K(x,E_2(x),\theta(x)), \, x\in U\setminus \mathcal{S}^+
$$ 
such that 
\begin{enumerate}
\item the angle $\angle (E_1(x),E_2(x))\ge\text{ const. }$;
\item $df(K^s(x))\subset K^s(f(x))$ for any $x\in U\setminus\mathcal{S}^+$ and
$df^{-1}(K^u(x))\subset K^u(f(x))$ for any $x\in f(U\setminus\mathcal{S}^+)$;
\item for some $\lambda>1$
\begin{enumerate}
\item $\|df_xv\|\ge\lambda\|v\|$ for $x\in U\setminus\mathcal{S}^+$ and $v\in K^u(x)$;
\item $\|df^{-1}_xv\|\ge\lambda\|v\|$ for $x\in f(U\setminus\mathcal{S}^+)$ and 
$v\in K^s(x)$.
\end{enumerate}
\end{enumerate}

\begin{theorem}[\cite{Pes2}]
Let $\Lambda$ be a hyperbolic attractor with singularities for a $C^{1+\alpha}$ map and assume that Condition \eqref{core} holds. Then $f$ admits an SRB measure on $\Lambda$. 
\end{theorem}

\subsection{Examples} We describe the following three examples of hyperbolic attractors with singularities which satisfy requirements \eqref{singularities} and \eqref{core} and thus possess SRB measures.

{\bf The Lorenz attractor}. Let $I=(-1,1)$, $U=I\times I$, $N=I\times {0}\subset U$ and $f:U\setminus N\to U$ is given by
$$
f(x,y)=((-B|y|^{\nu_0}+B\text{sign}(y)|y|^{\nu}+1)\text{sign}(y),\,((1+A)|y|^{\nu_0}-A)\text{sign}(y)),
$$
where 
$$
0<A<1, \quad 0<B<\frac12,\quad \nu>1, \quad \frac{1}{1+A}<\nu_0<1.
$$
This attractor appears in the Lorenz system of ODE :
$$
\dot{x}=-\sigma x+\sigma y, \quad \dot{y}=rx-y-xz,\quad \dot{z}=xy-bz
$$
for the values of the parameters $\sigma=10$, $b=\frac83$ and $r\sim 24.05$, see \cite{ABS, AP, BS, GW, Wil}.

{\bf The Lozi attractor}. Let $I=(-c,c)$ for some $0<c<1$ and let $U=I\times I$, 
$N={0}\times I\subset U$ and $f:U\setminus N\to U$ is given by 
$$
f(x,y)=(1+by-a|x|, x),
$$
where $0<a<a_0$ and $0<b<b_0$ for some small $a_0>0$ and $b_0>0$. 

Up to a change of coordinates this map was introduced by Lozi as a simple version of the famous H\'enon map in population dynamics, see \cite{Lev, Mis, Rych, You85}.

{\bf The Belykh attractor}. Let $I=(-1,1)$, $U=I\times I$, $N=\{(x,y):y=kx\}\subset U$ and 
$f:U\setminus N\to U$ is given by 
$$
f(x,y)=\begin{cases}
(\lambda_1(x-1)+1,\,\lambda_2(y-1)+1)&\text{ for } y>kx,\\
(\mu_1(x+1)-1,\,\mu_2(y+1)-1)&\text{ for } y<kx,
\end{cases}
$$
where
$$
0<\lambda_1, \mu_1<\frac12, \quad 1<\lambda_2, \mu_2<\frac{2}{1-|k|}, \quad |k|<1.
$$
In the case $\lambda_1=\mu_1$ and $\lambda_2=\mu_2$ this map was introduced by Belykh \cite{Bel} as one of the simplest models in the phase synchronization theory in radiophysics.

\end{document}